\documentclass[a4paper,11pt]{amsart}
\usepackage{amsmath,amsthm,amssymb,amsthm, amsfonts, amscd, mathrsfs,mathrsfs, xcolor}
\makeatletter
 \def\@makefnmark{%
 \leavevmode
 \raise.9ex\hbox{\check@mathfonts
 \fontsize\sf@size\z@\normalfont%
 \@thefnmark}%
 }
 \makeatother
\newcommand\diam{\textrm{diam}}

\newcommand{\N}{\mathbb{N}}
\newcommand{\R}{\mathbb{R}}
\newcommand{\SO}{\mathcal{SO}}

\theoremstyle{plain}
\newtheorem{thm}{Theorem}[section]
\newtheorem{lemma}[thm]{Lemma}

\newtheorem{corollary}[thm]{Corollary}
\theoremstyle{definition}

\newtheorem{example}[thm]{Example}
\newtheorem{claim}[]{Claim}
\newtheorem{remark}[thm]{Remark}
\newcounter{cn}
\title[A Banach--Stone-like Theorem for Slowly Oscillating Functions]{A Banach--Stone-like Theorem for Slowly Oscillating Functions}
\author{Yutaka Iwamoto}

\address{Faculty of Fundamental Science, National Institute of Technology (KOSEN), Niihama College,
Niihama, 792-8580, Japan}
\email{y.iwamoto@niihama-nct.ac.jp}
\begin{document}

\begin{abstract}
We prove that the homeomorphism induced by a lattice isomorphism between
the lattices of slowly oscillating real-valued functions on proper metric
spaces is a coarse homeomorphism.
This yields a Banach--Stone-like theorem for slowly oscillating functions
on arbitrary proper metric spaces, which removes the chain-connectedness
assumption in the author's previous work.
We also give a necessary and sufficient condition for a lattice isomorphism
to be continuous with respect to the topology of uniform convergence.
\end{abstract}

\maketitle

\section{Introduction}
The classical Banach--Stone theorem states that, for compact Hausdorff spaces
$X$ and $Y$, the Banach spaces $C(X)$ and $C(Y)$ of real-valued continuous
functions are linearly isometric if and only if $X$ and $Y$ are homeomorphic
\cite{Banach1,Stone}.
Kaplansky subsequently showed that $C(X)$ and $C(Y)$ are lattice-isomorphic
 if and only if $X$ and $Y$ are homeomorphic \cite{Kapl}.
Various extensions of these results have been obtained for lattices of
uniformly continuous functions
(see, e.g., \cite{FJ_Cabello_Sanchez,Garrido-Jaramillo-2000,Husek,H-P-1,H-P-2}).
For a survey of more general Banach--Stone theorems, see
\cite{GarridoJaramillo2002,LT}.

Coarse geometry studies the large-scale structure of metric spaces.
Let $X$ be a proper metric space.
A central object in this setting is the Higson compactification $hX$ and its remainder $\nu X=hX\setminus X$, called the Higson corona.
Denote by $\SO(X)$ the set of continuous real-valued slowly oscillating functions on $X$,
 and by $\SO_b(X)$ the subset of $\SO(X)$ consisting of all bounded functions.
The algebra $\SO_b(X)$ consists precisely of the bounded continuous functions
 on $X$ that extend continuously to $hX$ (see \cite{Keesling,Roe}; cf. \cite{Iwa2025}).

In \cite{Iwa2025}, the author established a Banach--Stone-like theorem
 for chain-connected proper metric spaces, showing that every lattice
 isomorphism from $\SO(X)$ onto $\SO(Y)$ induces a coarse homeomorphism between $X$ and $Y$.
However, the chain-connectedness assumption excludes many natural examples,
 for instance spaces that contain isolated points whose mutual distances tend to infinity.

A natural question is whether the chain-connectedness assumption can be removed entirely.
The present paper answers this question affirmatively,
 establishing both results for arbitrary proper metric spaces.

The main result is the following theorem.

\par\medskip
\noindent
\textbf{Banach--Stone-like theorem for slowly oscillating functions.}
(Theorem~\ref{thm_2-14})
\textit{Let $X$ and $Y$ be proper metric spaces.
The following statements are equivalent.}
\textit{
\begin{enumerate}
  \item $\SO(X)$ and $\SO(Y)$ are lattice-isomorphic.
  \item $\SO_b(X)$ and $\SO_b(Y)$ are lattice-isomorphic.
  \item There exists a coarse homeomorphism $\tau:X\to Y$.
  \item The Higson compactifications $hX$ and $hY$ are homeomorphic.
\end{enumerate}
}
\par\medskip

The key step in the proof is to show that the homeomorphism induced by a lattice isomorphism
 preserves direct summand sequences in the coarse category.
Here, we call a sequence of isolated points whose mutual distances tend to infinity
 a \textit{direct summand sequence in the coarse category}.
We show that the homeomorphism $\tau$ induced by a lattice isomorphism maps every such sequence in $X$ to such a sequence in $Y$.
Since $\tau^{-1}$ is induced by the inverse lattice isomorphism, the same holds for $\tau^{-1}$.
Consequently, both $\tau$ and $\tau^{-1}$ are coarse maps, that is, $\tau$ is a coarse homeomorphism,
 without any additional assumptions on $X$ and $Y$.

By analyzing these direct summand sequences, we also give a necessary
 and sufficient condition for a lattice isomorphism to be continuous in the topology of uniform convergence (Theorem~\ref{thm_3-01}).
As a consequence, we see that if neither $X$ nor $Y$ admits a direct summand sequence in the coarse category,
 then for every lattice isomorphism $T:  \SO(X)\to\SO(Y)$ with $T\underline{0}=\underline{0}$,
  the restriction of $T$ to $\SO_b(X)$ is a lattice isomorphism from $\SO_b(X)$ onto $\SO_b(Y)$ (Corollary~\ref{cor_3-02}).

\section{Banach--Stone-like Theorem}
Let $(X, d)$ be a metric space and let $B_r(x)$
 denote the closed ball of radius $r$ centered at $x\in X$.
A metric $d$ on $X$ is called \textit{proper} if $B_{r}(x)$ is compact
 for every $x\in X$ and $r>0$.
In this case, $(X,d)$ is called a \textit{proper metric space}.
Unless otherwise stated, metric spaces
$X$ and $Y$ are assumed to be unbounded proper metric spaces with base points $x_0$ and $y_0$, respectively.
To simplify notation, we write $|x|=d(x_0,x)$ for each $x\in X$.

A map $f:X\to Y$ from a noncompact proper metric space $X$ into a metric space $Y$
 is called a \textit{slowly oscillating function} if, for every $r>0$,
\[
  \lim_{|x|\to\infty}\textrm{diam}\bigl(f(B_r(x))\bigr)=0,
\]
  where $\diam A=\sup\{d(y,z)\mid y, z\in A\}$.
Equivalently, for every $r>0$ and every $\varepsilon>0$, there exists a compact subset $K\subset X$ such that
\[
  \textrm{diam}\bigl(f(B_r(x))\bigr)<\varepsilon
\]
  for every $x\in X\setminus K$.

Recall that $\SO(X)$ denotes the set of all continuous real-valued slowly oscillating functions on $X$,
 and that $\SO_b(X)$ denotes its subset of bounded functions.
For $f,g\in\SO(X)$, define the join $f\vee g$ and the meet $f \wedge g$ by
\[
(f\vee g)(x)=\max\{f(x),g(x)\},\quad
(f\wedge g)(x)=\min\{f(x),g(x)\}.
\]
Then both $\SO(X)$ and $\SO_b(X)$ are lattices with respect to these operations \cite{Iwa, Iwa2025}.

We use $\underline{a}$ to denote the constant function with value $a\in\R$.

\begin{example}\label{ex_2-1}
Let $A$ and $B$ be nonempty disjoint closed subsets of a proper metric space $(X,d)$ such that
 $d(A,B)=\inf\{d(x,y)\mid x\in A,~ y\in B\}>0$.
Define $\lambda_{A,B}:X\to[0,1]$ by
\[
\lambda_{A,B}(x)=\frac{d(x,B)}{d(x,A)+d(x,B)}.
\]
Then $\lambda_{A,B}$ is Lipschitz continuous, with $\lambda_{A,B}|_A=\underline{1}$ and $\lambda_{A,B}|_B=\underline{0}$.
In fact,
\[
|\lambda_{A,B}(x)-\lambda_{A,B}(y)|\leq \frac{1}{d(A,B)}d(x,y)
\]
for all $x,y\in X$.

Let $(x_i)$ be a sequence of points of $X$ and let $(r_i)$ be a sequence of positive numbers such that
\[
B_{r_i}(x_i)\cap B_{r_j}(x_j)=\emptyset, \qquad i\neq j.
\]
Put $A_n =\{ x_n\}$ and $B_n=\overline{X\setminus B_{r_n}(x_n)}$
 and let  $\lambda_n =\lambda_{A_n, B_n}$ for every $n$.
Consider the function 
\[
\lambda =\bigvee _{n=1}^{\infty} \lambda_n :X\to[0,1].
\]
\par\noindent
(i) If $r_n\to\infty$, then $\lambda$ is a bounded slowly oscillating function satisfying $\lambda(x_n)=1$ for every $n$
 and vanishing outside $\bigcup_{n=1}^{\infty}B_{r_n}(x_n)$.
\par\noindent
(ii) If $f\to\mathbb R$ is a bounded function that vanishes at infinity, then $f\cdot\lambda$ is a bounded slowly oscillating function satisfying $(f\cdot\lambda)(x_n)=f(x_n)$ for every $n$.
\hfill$\Box$
\end{example}

A function $\mathfrak{t} :X\times\mathbb{R}\to\mathbb{R}$ induces an operator $N_{\mathfrak{t}}$ defined by
\[
(N_{\mathfrak{t}} f)(x)=\mathfrak{t} (x,f(x))
\]
for every function $f : X\to \R$.
The operator $N_{\mathfrak{t}}$ is called the \textit{superposition (or Nemytskii) operator} \cite{AZ} induced by $\mathfrak{t}$.
A continuous map $\tau : X\to Y$ induces an operator $C_{\tau}: C(Y)\to C(X)$
 defined by
\[ (C_{\tau}g)(x)=g(\tau(x))=(g\circ \tau) (x) ,\qquad g\in C(Y),~ x\in X.
\]
The operator $C_{\tau}$ is called the \textit{composition operator} induced by $\tau$.

\begin{thm}\cite[Theorem 3.6]{Iwa2025}\label{thm_2-2}
  Let $X$ and $Y$ be proper metric spaces.
  If $T: \SO(X) \to \SO(Y)$ is a lattice isomorphism,
  then there exist a unique homeomorphism
  $\tau:X\to Y$ and a unique function
  $\mathfrak{t}:Y\times\mathbb{R}\to\mathbb{R}$ such that
  \begin{enumerate}
    \item $\mathfrak{t}(y,\cdot):\mathbb{R}\to\mathbb{R}$
    is an increasing homeomorphism for each $y\in Y$, and
    \item
    $Tf(y)=\mathfrak{t}(y,f(\tau^{-1}(y)))$
    for every $f\in\SO(X)$ and $y\in Y$.
  \end{enumerate}
  In other words, every lattice isomorphism between lattices
  of slowly oscillating functions is the composition of a superposition operator and a composition operator,
   namely,
\[
    T=N_{\mathfrak{t}}\circ C_{\tau^{-1}}.
\]
  The same conclusion holds for every lattice isomorphism
  $T:\SO_b(X)\to\SO_b(Y)$.
  \hfill$\Box$
\end{thm}

\begin{remark}\label{rem_2-5}
  Let $X$ and $Y$ be proper metric spaces.
  Since $\SO(X)\subset \R^{X}$, we can topologize it as a subspace of $\R^{X}$
   endowed with the Tychonoff product topology, called the topology of \textit{pointwise convergence}.
  A basic neighborhood of $f\in \R^{X}$ has the form
  \[
  N(f, F,\varepsilon )=\{g\in \R^{X}\mid |f(x)-g(x)|<\varepsilon \text{ for all } x\in F\},
  \]
  where $F$ is a finite subset of $X$ and $\varepsilon>0$.

  If $T :  \SO(X) \to \SO(Y)$ is a lattice isomorphism,
   then $T$ is a homeomorphism in the topology of pointwise convergence.
  Indeed, by Theorem~\ref{thm_2-2}, we have
  \[
  Tf(y)=\mathfrak{t}(y,f(\tau^{-1}(y)))
  \]
  for each $y\in Y$, where the function $\mathfrak{t}(y,\cdot) :  \R\to\R$ is an increasing homeomorphism.
  Hence $T$ is continuous in the topology of pointwise convergence.
  Since $T^{-1}$ is also a lattice isomorphism, the same argument shows that $T^{-1}$ is continuous.

  For the corresponding question for the topology of uniform convergence, see Theorem \ref{thm_3-01}.
  \hfill$\Box$
\end{remark}

Throughout this paper,
 the homeomorphism $\tau$ and the function $\mathfrak{t}$ appearing in the representation given by 
 Theorem~\ref{thm_2-2} are called the \textit{homeomorphism induced by $T$}
  and the \textit{representing function associated with $T$}, respectively.

\begin{corollary}\label{cor_2-4}
  Let $X$ and $Y$ be proper metric spaces.
  Let $T: \SO(X) \to \SO(Y)$ be a lattice isomorphism
  and let $\tau:X\to Y$ be the homeomorphism induced by $T$.
  Then $\tau^{-1}: Y\to X$ is the homeomorphism induced by $T^{-1}: \SO(Y)\to \SO(X)$.
  The same conclusion holds for every lattice isomorphism
  $T:\SO_b(X)\to\SO_b(Y)$.
\end{corollary}

\begin{proof}
  Let $\mathfrak{t}:Y\times\mathbb{R}\to\mathbb{R}$ be the representing
  function associated with $T$.
  For each $y\in Y$, we write $\mathfrak{t}_y=\mathfrak{t}(y,\cdot)$.
  By Theorem~\ref{thm_2-2}, $\mathfrak{t}_y:\mathbb{R}\to\mathbb{R}$
  is an increasing homeomorphism, and hence so is its inverse
  $\mathfrak{t}_y^{-1}$.

Since $T$ is a lattice isomorphism, so is $T^{-1}:\SO(Y)\to\SO(X)$,
and hence Theorem~\ref{thm_2-2} applies to $T^{-1}$.
  Define $\mathfrak{s}:X\times\mathbb{R}\to\mathbb{R}$ by
\[
    \mathfrak{s}(x,r)=\mathfrak{t}_{\tau(x)}^{-1}(r)
    \qquad (x\in X,\ r\in\mathbb{R}).
\]
  Then $\mathfrak{s}(x,\cdot)=\mathfrak{t}_{\tau(x)}^{-1}$ is an
  increasing homeomorphism of $\mathbb{R}$ for each $x\in X$.

  Now let $g\in\SO(Y)$ and put $f=T^{-1}g\in\SO(X)$.
  Then $g=Tf$, so by Theorem~\ref{thm_2-2},
\[
    g(y)=\mathfrak{t}_y\bigl(f(\tau^{-1}(y))\bigr)
    \qquad (y\in Y).
\]
  Substituting $y=\tau(x)$ for $x\in X$, we obtain
  $g(\tau(x))=\mathfrak{t}_{\tau(x)}\bigl(f(x)\bigr)$.
  Applying $\mathfrak{t}_{\tau(x)}^{-1}$ to both sides yields
\[
    T^{-1}g(x)=f(x)=\mathfrak{t}_{\tau(x)}^{-1}\bigl(g(\tau(x))\bigr)
    =\mathfrak{s}\bigl(x,g(\tau(x))\bigr)
    \qquad (x\in X).
\]
  That is, $T^{-1}=N_{\mathfrak{s}}\circ C_{\tau}$.
  Since $\tau=(\tau^{-1})^{-1}$, this shows that the homeomorphism
  $\tau^{-1}:Y\to X$ and the function $\mathfrak{s}$ satisfy
  conditions (1) and (2) of Theorem~\ref{thm_2-2} for the lattice
  isomorphism $T^{-1}:\SO(Y)\to\SO(X)$, that is,
  \begin{enumerate}
    \item $\mathfrak{s}(x,\cdot)$ is an increasing homeomorphism
    for each $x\in X$, and
    \item $T^{-1}g(x)=\mathfrak{s}\bigl(x,g((\tau^{-1})^{-1}(x))\bigr)$
    for every $g\in\SO(Y)$ and $x\in X$.
  \end{enumerate}
  By the uniqueness part of Theorem~\ref{thm_2-2}, $\tau^{-1}$ is the
  homeomorphism induced by $T^{-1}$ and $\mathfrak{s}$ is the
  representing function associated with $T^{-1}$.

  The same argument applies to lattice isomorphisms
  $T:\SO_b(X)\to\SO_b(Y)$.
\end{proof}

\begin{corollary}\label{cor_2-5}
  Let $X$ and $Y$ be proper metric spaces.
  Let $T: \SO(X) \to \SO(Y)$ be a lattice isomorphism
  and let $\tau:X\to Y$ be the homeomorphism induced by $T$.
  Then, for every $x\in X$, the following hold:
  \begin{enumerate}
    \item\label{cor_2-4-1}
    if $f,f'\in\SO(X)$ satisfy $f(x)=f'(x)$,
    then $Tf(\tau(x))=Tf'(\tau(x))$, and
    \item\label{cor_2-4-2}
    if $g,g'\in\SO(Y)$ satisfy $g(\tau(x))=g'(\tau(x))$,
    then $T^{-1}g(x)=T^{-1}g'(x)$.
  \end{enumerate}
  In other words, $Tf(\tau(x))$ depends only on $f(x)$, and
  $T^{-1}g(x)$ depends only on $g(\tau(x))$.
  The same conclusion holds for every lattice isomorphism
  $T:\SO_b(X)\to\SO_b(Y)$.
\end{corollary}

\begin{proof}
  Let $\mathfrak{t}$ be the representing function associated with $T$.
  If $f(x)=f'(x)$, then Theorem~\ref{thm_2-2} yields
\[
    Tf(\tau(x))=\mathfrak{t}\bigl(\tau(x),f(x)\bigr)
    =\mathfrak{t}\bigl(\tau(x),f'(x)\bigr)=Tf'(\tau(x)),
\]
  which proves (\ref{cor_2-4-1}).
  By Corollary~\ref{cor_2-4}, $\tau^{-1}$ is the homeomorphism induced by
  the lattice isomorphism $T^{-1}$.
  Hence (\ref{cor_2-4-2}) follows by applying (\ref{cor_2-4-1})
  to $T^{-1}$ in place of $T$.
\end{proof}

 For linear lattice isomorphisms,
  Theorem~\ref{thm_2-2} yields the following result.

\begin{thm}\cite[Theorem 3.8]{Iwa2025}\label{thm_2-6}
  Let $X$ and $Y$ be proper metric spaces.
  Let $T: \SO(X) \to \SO(Y)$ be a linear lattice isomorphism,
  and let $\tau :X \to Y$ be the homeomorphism induced by $T$.
  Then $T$ is a weighted composition operator of the form
  \[  T=\omega \cdot C_{\tau^{-1}},
  \]
   where $\omega=T(\underline{1})$.
  The same representation holds for every linear lattice isomorphism $T: \SO_b(X) \to \SO_b(Y)$.
  \hfill$\Box$
\end{thm}

Let $X$ be a proper metric space.
We say that $\N$ is a \textit{direct summand of $X$ in the coarse category} (cf. \cite{FJ_Cabello_Sanchez})
 if
\[
\sup \{ r\geq 0 \mid  B_r (x)=\{x\}\text{ for some } x\in X\}=\infty.
\]
If this condition holds,
 then there exist a sequence $(x_n)$ of points of $X$ and a sequence $(r_n)$ of positive numbers such that
\[
\lim_{n\to\infty} r_n=\infty
\quad\text{and}\quad
B_{r_n}(x_n)=\{x_n \}
\]
for every $n\in\mathbb{N}$.
We call any sequence $(x_n)$ of distinct poins satisfying these conditions a  \textit{direct summand sequence in $X$ in the coarse category}.
If $(x_n)$ is a direct summand sequence in $X$ in the coarse category then $|x_n|\to \infty$ as $n\to \infty$ since $X$ is a proper metric space.

\begin{example}\label{ex_2-7}
Let $X$ be a proper metric space.
Let $(x_n)$ be a direct summand sequence in $X$ in the coarse category.
For a slowly oscillating function $f:X\to\R$ and a sequence $(n_k)$ of real numbers,
 define a function $g: X\to\R$ by
\[
 g(x)=
 \begin{cases}
  n_k \cdot f(x_k)    & \textrm{if } x=x_k\\  
  f(x) & \text{if } x\not\in \{x_n\mid n\in\mathbb{N}\}.
 \end{cases}
\]
Then $g$ is a slowly oscillating continuous function.

The well-definedness and continuity of $g$
 follow from the fact that $|x_n|\to\infty$ and that each $x_n$ is isolated.
To see that $g$ is slowly oscillating, fix \(r>0\).
Let $(r_n)$ be a sequence of positive numbers such that
$\lim_{n\to\infty} r_n=\infty$ and $B_{r_n}(x_n)=\{x_n \}$.
There exists $n_0\in\mathbb{N}$
 such that 
 $r_n >r$ for every $n\geq n_0$.
Put
\[
 K=\bigcup_{i<n_0} B_{r_i} (x_i).
\]
Then $K$ is compact and
\[
\diam \, g(B_r(x))=0
\]
for every $x\in X\setminus K$.
Since $r>0$ was arbitrary, $g$ is slowly oscillating.
\hfill$\Box$
\end{example}

A metric space $X$ is said to be \textit{$r$-chain-connected} ($r>0$)
if, for any two points $x, x' \in X$,
there exist finitely many points $p_0, \dots, p_n$ in $X$ such that
$p_0 =x$, $p_n =x'$ and $d_X(p_{i-1}, p_{i})\leq r$ for every $i = 1,\dots, n$.
A maximal $r$-chain-connected subset of $X$ is called an \textit{$r$-chain-connected component} of $X$.
A metric space is said to be \textit{chain-connected} if it is $r$-chain-connected for some $r > 0$.

\begin{lemma}\label{lem_2-8}  
If $\mathbb{N}$ is not a direct summand of $X$ in the coarse category, then there exists $r>0$ such that every $r$-chain-connected component of $X$ contains at least two points.
\end{lemma}

\begin{proof}
  Since $\N$ is not a direct summand of $X$ in the coarse category, we have
\[
  M=\sup \{ R\geq 0 \mid  B_R (x)=\{x\},~ x\in X\}
\]
  is finite.
  Set $r=2M+1$.
  Then $B_r(x)$ contains at least two points for every $x\in X$.
  Hence, every $r$-chain-connected component of $X$ contains at least two points.
\end{proof}

A set is called \textit{non-degenerate} if it contains at least two points.

\begin{corollary}\label{cor_2-9}
  If a sequence $(x_n)$ is not a direct summand sequence in $X$ in the coarse category,
   then there exists $r>0$ such that every $x_n$ belongs to a non-degenerate $r$-chain-connected component of $X$.
  \hfill$\Box$
\end{corollary}

The following lemma is fundamental to our study.
It shows that the homeomorphism induced by a lattice isomorphism preserves
 direct summand sequences in the coarse category of proper metric spaces.

\begin{lemma}\label{lem_2-10} 
  Let $X$ and $Y$ be proper metric spaces.
  Let $T: \SO(X) \to \SO(Y)$ be a lattice isomorphism with $T\underline{0}=\underline{0}$,
   and let $\tau : X\to Y$ be the homeomorphism induced by $T$.
  Then a sequence $(x_n)$ in $X$ is a direct summand sequence in the coarse category
   if and only if $(\tau(x_n))$ is a direct summand sequence in $Y$ in the coarse category.
\end{lemma}

\begin{proof}
  We prove only the forward implication.
  The reverse implication follows by applying the same argument to $T^{-1}$ and its induced homeomorphism using Corollary~\ref{cor_2-4}.

  Let $(x_n)$ be a direct summand sequence in $X$ in the coarse category.
  Suppose, on the contrary, that $(\tau(x_n))$ is not a direct summand sequence in $Y$ in the coarse category.
  Put $y_n =\tau(x_n)$ for each $n\in \N$.
Then $|y_n|\to \infty$ as $n\to \infty$, since $|x_n|\to \infty$ and $\tau$ is a homeomorphism between proper metric spaces.

  Let $f=T^{-1}\underbar{1}\in \SO(X)$.
  Then we define $h:X\to \R$ by $h(x_n)=f(x_n)$ and $h(x)=0$ whenever $x\not \in \{x_n \mid n\in\N\}$.
  Since $(x_n)$ is a direct summand sequence in $X$ in the coarse category,
   $h$ is a well-defined slowly oscillating continuous function (see Example~\ref{ex_2-7}).
  Put $g=Th\in \SO(Y)$.
Since $h(x_n)=f(x_n)$, Corollary~\ref{cor_2-5} yields that
  \begin{align*}
    g(y_n)
    &=Th(y_n)=Tf(y_n)=1.
  \end{align*}
  On the other hand, if $y=\tau (x)\notin \{y_n\mid n\in\N\}$, then
  $h(x)=0$ implies that
  \begin{align*}
    g(y)
    &=Th(y)=T\underbar{0}(y)=0.
  \end{align*}
Since $(y_n)$ is not a direct summand sequence in $Y$ in the coarse category,
  by Corollary~\ref{cor_2-9},
   there exists $r>0$ such that each $y_n $ belongs to a non-degenerate $r$-chain-connected component of $Y$.
After passing to a subsequence of $(y_n)$ if necessary,
 we can take a sequence $(z_n)$ in $Y$ such that $0<d(y_n, z_n)\leq r$
 with $z_n \not\in \{y_k \mid k\in\N\}$.
 Consequently, $g(y_n)=1 $ and $g(z_n)=0$ for all $n$,
  contradicting the assumption that $g\in\SO(Y)$.
\end{proof}

Let $f: X\to Y$ be a map between metric spaces $(X,d_X)$ and $(Y,d_Y)$.
 We call $f$ \textit{uniformly expansive} (or \textit{uniformly bornologous})
  if there exists a non-decreasing function $\sigma:[0,\infty)\to[0,\infty)$ such that
\[
d_Y(f(x),f(x'))\leq \sigma(d_X(x,x'))
\]
for all $x,x'\in X$, and \textit{metrically proper}
 if the preimage of every $d_Y$-bounded subset of $Y$ is $d_X$-bounded in $X$.
Every homeomorphism between proper metric spaces is metrically proper.

A map that is both uniformly expansive and metrically proper is called a \textit{coarse map}. 
Maps $f,g: X\to Y$ are \textit{close} if 
\[
\sup_{x\in X } d_Y(f(x),g(x)) <\infty. 
\]
A coarse map $f: X\to Y$ is a \textit{coarse equivalence} if there exists a coarse map $g: Y\to X$
 such that $g\circ f$ and $f\circ g$ are close to $\mathrm{id}_X$ and $\mathrm{id}_Y$, respectively. 
In this case, $X$ and $Y$ are said to be \textit{coarsely equivalent}. 
Coarsely equivalent proper metric spaces have homeomorphic Higson coronas \cite[Corollary 2.42]{Roe}.

A homeomorphism $f\to Y$ is called a \textit{coarse homeomorphism} if both $f$ and $f^{-1}$ are coarse maps.

The following is a fundamental relation between coarse maps and slowly oscillating functions.

\begin{lemma}\cite[Lemma 4.1]{Iwa2025}\label{lem_2-11}
  If $h:X \to Y$ is a coarse continuous map between proper metric spaces
  and $f: Y\to \R$ is slowly oscillating
  then $f\circ h: X\to \R$ is a slowly oscillating function.
  \hfill$\Box$
\end{lemma}

The following lemma removes the chain-connectedness assumption from \cite[Lemma 4.2]{Iwa2025}.
It is a key ingredient in establishing the main results of \cite{Iwa2025} without this assumption.

\begin{lemma}\label{lem_2-12} 
  Let $X$ and $Y$ be proper metric spaces
  and let $T: \SO(X) \to \SO(Y)$ be a lattice isomorphism with $T(\underline{0}) = \underline{0}$.
  Then the homeomorphism $\tau: X \to Y$ induced by $T$ is uniformly expansive.
  The same conclusion holds for lattice isomorphisms $T: \SO_b(X) \to \SO_b(Y)$
   with $T\underline{0}=\underline{0}$. 
\end{lemma}

\begin{proof}
  We prove the lemma only for $T: \SO(X) \to \SO(Y)$.
  The same argument applies to lattice isomorphisms
  $T: \SO_b(X) \to \SO_b(Y)$.

  Note also that $T^{-1}\underline{0}=\underline{0}$
  since $T\underline{0}=\underline{0}$.

  Suppose, to the contrary, that $\tau$ is not uniformly expansive.
  Then there exist a constant $R>0$ and sequences $(a_n)$ and $(b_n)$
  in $X$ such that
  \begin{enumerate}
    \item\label{L-1} $0 < d(a_n, b_n) \leq R$ for every $n$, and
    \item\label{L-2} $d(\tau(a_n), \tau(b_n)) \to \infty$ ($n \to \infty$).
  \end{enumerate}
  Since $X$ is a proper metric space and $\tau$ is metrically proper,
  (\ref{L-1}) and (\ref{L-2}) imply that
  \begin{enumerate}\setcounter{enumi}{2}
    \item\label{L-3} $|a_n| \to \infty$ and $|b_n| \to \infty$, and
    \item\label{L-4} $|\tau(a_n)| \to \infty$ and $|\tau(b_n)| \to \infty$
    ($n \to \infty$).
  \end{enumerate}
  Passing to subsequences, we may assume that
  \begin{enumerate}\setcounter{enumi}{4}
    \item\label{L-5} $|\tau(a_n)| + 2n^2 < |\tau(b_n)| < |\tau(a_{n+1})| - 2(n+1)^2$
    for every $n$.
  \end{enumerate}
  Hence we can take $g \in \SO_b(Y)$ such that
  $g(\tau(a_n)) = 0$ and $g(\tau(b_n)) = 1$ for every $n$
  (cf.\ Example~\ref{ex_2-1}).

  Put $f = T^{-1}g \in \SO(X)$.
  Since $g(\tau(a_n)) = 0 = \underline{0}(\tau(a_n))$,
  Corollary~\ref{cor_2-5} yields
  $f(a_n) = T^{-1}\underline{0}(a_n) = 0$ for every $n$.
  Since $f$ is slowly oscillating, it follows from (\ref{L-1}) and
  (\ref{L-3}) that $f(b_n) \to 0$ ($n \to \infty$).

  By (\ref{L-1}), the sequence $(b_n)$ is contained in the union of all
  non-degenerate $R$-chain-connected components of $X$.
  By Lemma~\ref{lem_2-10}, $(\tau(b_n))$ is not a direct summand sequence
  in $Y$ in the coarse category.
  Thus, Corollary~\ref{cor_2-9} yields $r > 0$ and a sequence
  $(c_n)$ in $X$ such that
  \begin{enumerate}\setcounter{enumi}{5}
    \item\label{L-6} $0 < d(\tau(b_n), \tau(c_n)) < r$ for every $n$.
  \end{enumerate}
  By (\ref{L-5}), we may also assume
  that $\tau(c_n) \neq \tau(b_m)$ for every $m, n$.
  Since $f(b_n) \to 0$ ($n \to \infty$), there exists $h \in \SO_b(X)$
  such that
\[
    h(c_n) = 0 \quad\text{and}\quad h(b_n) = f(b_n)
    \quad\text{for every } n.
\]
  By Corollary~\ref{cor_2-5}, we obtain
\[
    Th(\tau(c_n)) = T\underline{0}(\tau(c_n)) = 0
    \quad\text{and}\quad
    Th(\tau(b_n)) = Tf(\tau(b_n)) = g(\tau(b_n)) = 1
\]
  for every $n$.
  On the other hand, since $Th \in \SO(Y)$, it follows from
  (\ref{L-4}) and (\ref{L-6}) that
  $|Th(\tau(c_n)) - Th(\tau(b_n))| \to 0$ ($n \to \infty$),
  a contradiction.
\end{proof}

\begin{corollary}\label{cor_2-13}
  Let $X$ and $Y$ be proper metric spaces
  and let $T: \SO(X) \to \SO(Y)$ be a lattice isomorphism with $T\underline{0} = \underline{0}$.
  Then the homeomorphism $\tau :X \to Y$ induced by $T$ is a coarse homeomorphism.
The same conclusion holds for lattice isomorphisms $T:\SO_b(X)\to\SO_b(Y)$ satisfying $T\underline{0}=\underline{0}$.
\end{corollary}

\begin{proof}
  Since a homeomorphism between proper metric spaces is metrically proper,
  $\tau$ is a coarse homeomorphism by Lemma~\ref{lem_2-12}.
\end{proof}

Let 
\[
e_X : X\to \displaystyle\prod_{f\in\SO_b(X)}[\inf f,\sup f]
\]
 be the evaluation map of $X$ associated with $\SO_b(X)$,
 i.e., 
\[
e_X (x)=(f(x))_{f\in\SO_b(X)}.
\]
This map is a topological embedding, and its closure $\overline{e_X (X)}$ is a compactification of $X$. 
This compactification $hX=\overline{e_X (X)}$ and its remainder $\nu X =hX\setminus X$ are called the \textit{Higson compactification} and the \textit{Higson corona} of $X$, respectively. 
In coarse geometry, these constructions are used to study the geometric behavior of proper metric spaces near infinity \cite{Roe}.
For basic properties of Higson compactifications, see \cite{Keesling, Iwa2025}.

The preceding results show that all the results used to prove the main theorem of \cite[Theorem 4.5]{Iwa2025} remain valid
 without the chain-connectedness assumption.
We therefore obtain the following theorem for arbitrary proper metric spaces.

\begin{thm}[Banach--Stone-like theorem for slowly oscillating functions]\label{thm_2-14}
  Let $X$ and $Y$ be proper metric spaces.
  The following are equivalent.
  \begin{enumerate}
    \item\label{BS_1} $\SO(X)$  and $\SO(Y)$ are lattice isomorphic.
    \item\label{BS_2} $\SO_b(X)$  and $\SO_b(Y)$ are lattice isomorphic.
    \item\label{BS_3} There exists a coarse homeomorphism between $X$ and $Y$.
    \item\label{BS_4} The Higson compactifications $hX$ and $hY$ are homeomorphic.
  \end{enumerate}
\end{thm}

\begin{proof}
We shall show that $(1)\Leftrightarrow(3)$, $(2)\Rightarrow(3)$, $(3)\Rightarrow(4)$ and $(4)\Rightarrow(2)$.

  The implications $(\ref{BS_1}) \Rightarrow (\ref{BS_3})$ and $(\ref{BS_2}) \Rightarrow (\ref{BS_3})$ follow from Corollary~\ref{cor_2-13}.
  Conversely, by Lemma~\ref{lem_2-11},
   every coarse homeomorphism $\tau : X \to Y$ 
   induces a lattice isomorphism
  $C_{\tau}|_{\SO(Y)}: \SO(Y) \to \SO(X).$
  Thus, $(\ref{BS_3}) \Rightarrow (\ref{BS_1})$.
  
  $(\ref{BS_3}) \Rightarrow (\ref{BS_4})$.
  Let $\tau: X \to Y$ be a coarse homeomorphism.
  Let $e_Y : Y \to hY$ be the evaluation map associated with $\SO_b(Y)$.
  By Lemma~\ref{lem_2-11}, $f\circ\tau\in\SO_b(X)$ for every $f\in\SO_b(Y)$.
  Hence, by the characterization of Higson compactifications,
  the map $e_Y\circ\tau:X\to hY$ extends uniquely to a continuous
  map $\hat{\tau}:hX\to hY$.
  Similarly, $\tau^{-1}:Y\to X$ extends to a continuous map
  $\widehat{\tau^{-1}}:hY\to hX$.
 Since $(\widehat{\tau^{-1}}\circ\hat{\tau})|_X=\textrm{id}_{X}$
   and $(\hat{\tau}\circ\widehat{\tau^{-1}})|_Y=\textrm{id}_{Y}$,
   and since $X$ and $Y$ are dense in $hX$ and $hY$, respectively,
   we have
  \[
    \widehat{\tau^{-1}}\circ\hat{\tau}=\textrm{id}_{hX},
    \qquad
    \hat{\tau}\circ\widehat{\tau^{-1}}=\textrm{id}_{hY}.
  \]
  Therefore, $\hat{\tau}:hX\to hY$ is a homeomorphism, i.e., $(\ref{BS_3}) \Rightarrow (\ref{BS_4})$.
\[
    \begin{CD}
      hX @>{\hat{\tau}}>> hY \\
      @A{e_X }AA   @AA{e_Y}A \\
      X   @>{\tau}>> Y
    \end{CD}
\]

 $(\ref{BS_4}) \Rightarrow (\ref{BS_2})$.
  Suppose that $(\ref{BS_4})$ holds, and let $\hat{\tau}:hX\to hY$ be a homeomorphism.
  By \cite[Proposition 4.1]{LT} (cf. \cite[Lemma 1]{Garrido-Jaramillo-2000}), we have $\hat{\tau}(X) = Y$.
  Thus, $\tau = \hat{\tau}|_{X}$ is a homeomorphism from $X$ onto $Y$.
For each $f\in\SO_b(Y)$, let $\hat{f}:hY\to\mathbb{R}$
  be its continuous extension.
  Then $f\circ\tau$ extends continuously to $hX$ via
  $\hat{f}\circ\hat{\tau}$.
  Hence, $f\circ\tau\in\SO_b(X)$ by the characterization of
  Higson compactifications.
  Thus, $C_{\tau}$ maps $\SO_b (Y)$ into $\SO_b (X)$.
  Applying the same argument to $\tau^{-1}$ shows that
  $C_\tau|_{\SO_b(Y)}:\SO_b(Y)\to\SO_b(X)$ is a lattice
  isomorphism. Therefore,
  $(\ref{BS_4})\Rightarrow(\ref{BS_2})$.
\end{proof}

\section{Lattice Isomorphisms Preserving Uniform Convergence}
In this section, we characterize the proper metric spaces for which every lattice isomorphism between the lattices of slowly oscillating functions is a homeomorphism in the topology of uniform convergence.
In general, a lattice isomorphism $T :  \SO(X) \to \SO(Y)$ need not be continuous in this sense, and the obstruction turns out to be precisely the presence of $\N$ as a direct summand in the coarse category.
Moreover, this continuity is equivalent to the condition that $T$ preserves bounded differences of functions (Theorem~\ref{thm_3-01}).
The proof of Theorem~\ref{thm_3-01} is obtained by modifying the proof of \cite[Corollary 2]{FJ_Cabello_Sanchez}.
As an application, it follows that, under the same assumption on $X$ and $Y$, every lattice isomorphism $T$ with $T\underline{0}=\underline{0}$ restricts to a lattice isomorphism between the lattices of bounded slowly oscillating functions (Corollary~\ref{cor_3-02}).

\begin{thm}\label{thm_3-01}
Let $X$ and $Y$ be proper metric spaces such that $\SO(X)$ and $\SO(Y)$ are lattice isomorphic. Then the following are equivalent.
\begin{enumerate}
\item $\mathbb{N}$ is not a direct summand of either $X$ or $Y$ in the coarse category.
\item Every lattice isomorphism $T:\SO(X)\to\SO(Y)$ preserves bounded differences, i.e., for all $f,g\in\SO(X)$,
\[
g-f\text{ is bounded}
\quad\Longleftrightarrow\quad
Tg-Tf\text{ is bounded}.
\]
\item Every lattice isomorphism $T:\SO(X)\to\SO(Y)$ is a homeomorphism in the topology of uniform convergence, i.e., for all $f\in\SO(X)$ and all sequences $(f_n)$ in $\SO(X)$,
\[
f_n\to f\text{ uniformly}
\quad\Longleftrightarrow\quad
Tf_n\to Tf\text{ uniformly}.
\]
\end{enumerate}
\end{thm}

\begin{proof}[Proof of $(2)\Rightarrow(1)$ and $(3)\Rightarrow(1)$ in Theorem~\ref{thm_3-01}]
  Suppose that $\N$ is a direct summand of $X$ or $Y$ in the coarse category.
  Since conditions (2) and (3) are invariant under replacing $T$ with $T^{-1}$, we may assume that $\N$ is a direct summand of $X$.
  Let $(x_n)$ be a direct summand sequence in $X$ in the coarse category.
  Let $T:\SO(X)\to \SO(Y)$ be a lattice isomorphism.
  We may assume that $T\underline{0}=\underline{0}$.
  Let $\tau : X\to Y$ be the induced homeomorphism associated with $T$.
  Put $y_n =\tau (x_n)$ for each $n\in \N$.
  Then $(y_n)$ is a direct summand sequence in $Y$ in the coarse category by Lemma~\ref{lem_2-10}.
For each $f\in \SO(X)$, define $Sf$ by
\[
 Sf(y)=
 \begin{cases}
  n\cdot Tf(y_n) & \text{if } y=y_n,\\
  Tf(y)    & \textrm{if } y\not\in \{y_n \mid n\in \N \} .
 \end{cases}
\]
As we have seen in Example~\ref{ex_2-7},
 $Sf$ is a well-defined slowly oscillating function for every $f\in \SO(X)$.
Then $S$ maps $\SO(X)$ into $\SO(Y)$.
Now it is easy to see that $S$ is a lattice isomorphism such that $S\underline{0}=\underline{0}$.

To see $(2)\Rightarrow(1)$, 
 put $f=T^{-1}\underline{1}$.
Then $T^{-1}\underline{1} -T^{-1}\underline{0}=f-\underline{0} =f$ is bounded by our assumption.
However $Sf-S\underline{0}=Sf$ is unbounded
 since $Sf(y_n)=n\cdot Tf(y_n)=n\cdot 1 =n$ for every $n\in \N$.
This completes the proof of $(2)\Rightarrow(1)$.

To see $(3)\Rightarrow(1)$,
 define $g_n : Y\to \R$ by
\[
 g_n(y)=
 \begin{cases}
  \frac1{n} & \text{if } y=y_n,\\
  0  & \textrm{if } y\not\in \{y_n \mid n\in \N \} .
 \end{cases}
\]
Then $g_n$ is a slowly oscillating function (cf. Example~\ref{ex_2-7})
 and $g_n \to \underline{0}$ uniformly.
Put $f_n =T^{-1}g_n$ for each $n\in \N$.
Then $f_n \to T^{-1}\underline{0}=\underline{0}$ uniformly by our assumption.
However, $ (Sf_n)$ cannot converge to $S\underline{0}=\underline{0}$ uniformly since 
\[
Sf_n (y_n)= n\cdot Tf_n (y_n) = n\cdot g_n (y_n)=n\cdot \frac1{n}=1
\]
 for every $n\in \N$.
This completes the proof of $(3)\Rightarrow(1)$.
\end{proof}

\begin{proof}[Proof of $(1)\Rightarrow(2)$ in Theorem~\ref{thm_3-01}]
  Let $T : \SO(X)\to \SO(Y)$ be a lattice isomorphism 
  and let
  $\tau :X \to Y$ be the induced homeomorphism.
  Suppose that there exist $f, g\in \SO(X)$ such that $g-f$ is bounded but $Tg-Tf$ is unbounded.
  The map $h\mapsto Th-T\underline{0}$ is again a lattice isomorphism and $(Tg-T\underline{0})-(Tf-T\underline{0})=Tg-Tf$ is unbounded,
  we may assume that $T\underline{0}=\underline{0}$.

Replacing $f$ and $g$ with $f\wedge g$ and $f\vee g$, respectively,
 we may assume that $g-f>0 $.
Put 
$w=g-f$.
Then $w\in \SO(X)$ and $0<w<\infty$.

\begin{claim}\label{claim_02}
For $s\in[0,1]$ put $f_s=f+s\,w=(1-s)f+s\,g$. Then $f_s\in \SO(X)$, $f_0=f$, $f_1=g$, and
\[
  Tf\le Tf_s\le Tf_{s'}\qquad (0\le s\le s'\le 1).
\]
Consequently the set
\[
  I=\{\,s\in[0,1] : Tf_s-Tf \text{ is bounded}\,\}
\]
is an interval containing $0$, and $1\notin I$. 
\end{claim}
\smallskip
\noindent\textit{Proof.}
Since $w=g-f\in \SO(X)$, we have $f_s\in \SO(X)$. If $s\le s'$, then $f_s\le f_{s'}$ because $w >0$.
Since $T$ preserves the order,
 we have $Tf\le Tf_s\le Tf_{s'}$ for $0\le s\le s'\le 1$.
If $s\in I$ and $0\le s'\le s$, then $0\le Tf_{s'}-Tf\le Tf_s-Tf$, so $s'\in I$. Hence $I$ is an interval containing $0$.
Finally $f_1=g$ and $Tg-Tf$ is unbounded, so $1\notin I$. \hfill$\diamondsuit$

\medskip

Put $a=\sup I$.

\medskip
\begin{claim}
  There exist $\sigma\in\{+1,-1\}$ and a strictly decreasing sequence $\varepsilon_n \downarrow 0$
 of positive numbers such that $s_n=a+\sigma\varepsilon_n\in[0,1]$ and the functions
\[
  \varphi_n=\sigma\,(Tf_{s_n}-Tf_a)
\]
are unbounded above with
\[
  0\le \varphi_{n+1}\le\varphi_n\le\varphi_1
\]
for every $n\in\mathbb N$.

\end{claim}

\smallskip
\noindent\textit{Proof.}
\emph{Case 1: $a\in I$.} Then $a<1$ since $1\notin I$, and $I=[0,a]$. Put $\sigma=+1$
 and choose a strictly decreasing sequence
  $\varepsilon_n \downarrow 0$ with $a+\varepsilon_n\le 1$.
  Then $s_n=a+\varepsilon_n\notin I$ and $s_n\downarrow a$.
Thus $Tf_{s_n}-Tf$ is unbounded while $Tf_a-Tf$ is bounded, so
\[
\varphi_n=Tf_{s_n}-Tf_a
\]
 is unbounded above. By Claim~\ref{claim_02}, $\varphi_n\ge 0$, and $\varphi_n$ decreases in $n$ because $s_n$ decreases.

\emph{Case 2: $a\notin I$.} Then $I=[0,a)$ and $a>0$ (as $0\in I$). Put $\sigma=-1$
 and choose a strictly decreasing sequence
  $\varepsilon_n \downarrow 0$ with $\varepsilon_n<a$.
  Then $s_n=a-\varepsilon_n\in I$ and $s_n\uparrow a$.
Now $Tf_a-Tf$ is unbounded while $Tf_{s_n}-Tf$ is bounded, so
\[
\varphi_n=Tf_a-Tf_{s_n}
\]
 is unbounded above. 
By Claim~\ref{claim_02}, $\varphi_n\ge0$, and $\varphi_n$ decreases in $n$ because $s_n$ increases. \hfill$\diamondsuit$

\medskip
\begin{claim}
There exists a sequence $(y_n)$ in $Y$ such that $|y_n|\to\infty$ and $\varphi_n(y_n)\ge n$ for every $n$.  
\end{claim}

\smallskip
\noindent\textit{Proof.}
Since $Y$ is proper, $B_n(y_0)$ is compact.
Put $c_n=\sup_{B_n(y_0)}\varphi_1$. Then $0\le c_n<\infty$.
Since $\varphi_n$ is unbounded above, we can choose $y_n\in Y$ with $\varphi_n(y_n)\ge n+c_n$.
If $y_n\in B_n(y_0)$, then $\varphi_n(y_n)\le\varphi_1(y_n)\le c_n$, a contradiction.
Hence $|y_n|>n$, and $\varphi_n(y_n)\ge n$. \hfill$\diamondsuit$

\medskip
Since $\mathbb N$ is not a direct summand of $Y$ in the coarse category,
 there exists $r>0$ such that $B_r(y)\ne\{y\}$ for every $y\in Y$ (Lemma~\ref{lem_2-8} ).
Pick $z_n\in Y$ with $0< d(y_n, z_n)\le r$.

Put $a_n =\tau^{-1}(y_n)$ and $b_n=\tau^{-1}(z_n)$ for every $n\in \N$.
Note that $|a_n|\to \infty$ since $|y_n|\to\infty$.
Since $\tau^{-1}$ is uniformly expansive (Lemma~\ref{lem_2-12}),
 there exists $R>0$ such that $0< d(a_n , b_n)<R$ for every $n\in \N$.
Put $\delta_n =d(a_n, b_n)$ for every $n\in \N$.
Passing to a subsequence of $(a_n)$, we may assume
\[
  |a_{n+1}|\ge |a_n|+2R+1 \qquad (n\in\mathbb N),
\]
so that $d(a_n,a_m)>2R$ for $n\ne m$. 
In particular the balls $B_R(a_n)$ are pairwise disjoint.
Define $\lambda_n :  X\to[0,1]$ by
\[
  \lambda_n(x)=\max\bigl\{0,\;1-d(x,a_n)/\delta_n\bigr\}.
\]
Then $\lambda_n$ is continuous, $\lambda_n(a_n)=1$, $\lambda_n(b_n)=0$, and $\lambda_n=0$ outside $B_{\delta_n}(a_n)\subset B_R(a_n)$.
Put
\[
  \eta'=\sum_{n=1}^\infty \varepsilon_n\lambda_n .
\]
The supports are pairwise disjoint, and every compact set meets only finitely many of the balls $B_R(a_n)$
 (because $|a_n|\to\infty$ and $X$ is proper), so $\eta'$ is continuous.
Moreover, $0\le\eta'\le\varepsilon_1$, $\eta'(a_n)=\varepsilon_n$ and $\eta'(b_n)=0$.
Since $\varepsilon_n \downarrow 0$, it follows that $\eta' \in \SO_b (X)$.

Now put
\[
  \eta=f_a+\sigma\,\eta'\,w .
\]
Both $\eta'$ and $w$ are bounded slowly oscillating functions,
so $\eta'w\in \SO_b(X)$. Hence $\eta\in \SO(X)$.
Since $s_n=a+\sigma\varepsilon_n$, we have
\begin{align*}
  \eta(a_n)
  &=f_a (a_n) +\sigma \varepsilon_n w(a_n)
  =(f (a_n) +aw(a_n))+(s_n -a)w(a_n)\\
  &=f(a_n)+s_n w(a_n) =f_{s_n}(a_n)
\end{align*}
and $\eta(b_n)=f_a(b_n)$.
By Corollary~\ref{cor_2-5}, $\eta(a_n)=f_{s_n}(a_n)$ implies that
 $
  T\eta(y_n)
  = Tf_{s_n}(y_n)
 $.
Similarly, we have $T\eta(z_n)=Tf_a(z_n)$.
Therefore
\begin{align*}
  T\eta(y_n)-T\eta(z_n)
  &=Tf_{s_n}(y_n) - Tf_a (z_n)\\
  &=\left(Tf_{s_n}(y_n)-Tf_a (y_n) \right)+Tf_a (y_n)-Tf_a (z_n)\\
  &=\sigma \, \varphi_n (y_n) +\left( Tf_a (y_n)-Tf_a (z_n) \right).
\end{align*}
Because $Tf_a\in \SO(Y)$, $|y_n|\to\infty$ and $d(y_n, z_n)\le r$, $Tf_a (y_n)-Tf_a (z_n) $ tends to $0$,
 whereas $|\varphi_n(y_n)|\ge n\to\infty$.
Thus $|T\eta(y_n)-T\eta(z_n)|\to\infty$.
On the other hand, $T\eta\in \SO(Y)$ gives $|T\eta(y_n)-T\eta(z_n)|\to0$, a contradiction.
Hence $Tg-Tf$ is bounded whenever $g-f$ is bounded.
Applying the same argument to $T^{-1}:\SO(Y)\to\SO(X)$, and using that $\mathbb{N}$ is not a direct summand of $X$ in the coarse category, we obtain the converse implication.
This completes the proof of $(1)\Rightarrow(2)$.
\end{proof}

\medskip
\begin{proof}[Proof of $(1)\Rightarrow(3)$ in Theorem~\ref{thm_3-01}]
  Suppose, on the contrary, that there exists a lattice isomorphism
 $T:\SO(X)\to\SO(X)$ that is not continuous in the topology of uniform convergence.
Let
  $\tau :X \to Y$ be the induced homeomorphism.

The map $h\mapsto Th-T\underline{0}$ is again a lattice isomorphism and is continuous in the topology of uniform convergence if and only if $T$ is.
Hence we may assume $T\underline{0}=\underline{0}$.

Take a sequence $(f_n)\subset \SO(X)$ and a function $f\in \SO(X)$
 such that $(f_n)$ converges uniformly to $f$, but $(Tf_n)$ does not converge uniformly to $Tf$.
Put $f'_n=f_n-f$.
Then $\|f'_n\|_\infty\to0$.
Define $T'h=T(h+f)-Tf$ for each $h\in \SO(X)$.
Then $T': \SO(X)\to \SO(Y)$ is again a lattice isomorphism such that
 $T'\underline{0}=\underline{0}$, and $T'f'_n=Tf_n-Tf$ does not converge uniformly to $\underline{0}$.
Passing to a subsequence, we may assume that there are $\varepsilon_0>0$ and points $y_n\in Y$ such that
\[
  |T'f'_n(y_n)|>\varepsilon_0\qquad(n\in\mathbb N).
\]
Passing to a further subsequence, we may assume that either $T'f'_n(y_n)>\varepsilon_0$ for all $n$, or $T'f'_n(y_n)<-\varepsilon_0$ for all $n$.
In the second case we replace $T'$ by $S :  h\mapsto -T'(-h)$, and $f'_n$ by $-f'_n$.
Then $S$ is again a lattice isomorphism such that $S\underline{0}=\underline{0}$, and $Sf'_n(y_n)>\varepsilon_0$.
Relabeling, we have shown that we may assume the following:
\[
   T\underline{0}=\underline{0},\quad \|f'_n\|_\infty\to0,\quad Tf'_n(y_n)>\varepsilon_0\ \ (n\in\mathbb N).
\]
Now put
\[
  c_n=\sup_{k\ge n}\|f'_k\|_\infty+\tfrac1n 
\]
for every $n\in \N$.
Then $c_n \downarrow 0$ and $f'_n\le \underline{c_n}$.
Since $T$ preserves the order,
\[
  T\underline{c_n}(y_n)\ \ge\ Tf'_n(y_n)\ >\ \varepsilon_0\qquad(n\in\mathbb N).
\]
Note that $T\underline{c_n}(y)$ is non-negative and decreasing in $n$ for each $y\in Y$.

\medskip
\begin{claim}
   $|y_n|\to\infty$.
\end{claim}
\smallskip
\noindent\textit{Proof.}
Suppose that $(y_n)$ has a subsequence $(y_{n_k})$ converging to some $y\in Y$.
Fix $n\in\mathbb N$.
For $n_k\ge n$, we have $c_{n_k}\le c_n$, and hence
\[
  T\underline{c_n}(y_{n_k})\ge  T\underline{c_{n_k}}(y_{n_k})>\varepsilon_0 .
\]
Since $T\underline{c_n}\in \SO(Y)$ is continuous, letting $k\to\infty$ gives $T\underline{c_n}(y)\ge\varepsilon_0$ for every $n$.
However, by Remark~\ref{rem_2-5} and Corollary~\ref{cor_2-5},
 $c_n \downarrow 0$ implies that $T\underline{c_n}(y)\to T\underline{0}(y)=0$ as $n\to\infty$.
This is a contradiction.
Thus $(y_n)$ has no convergent subsequence. Since $Y$ is proper, every bounded sequence has a convergent subsequence, and hence $|y_n|\to\infty$.
\hfill$\diamondsuit$

\medskip
By passing to a subsequence of $(y_n)$ if necessary, and since $\mathbb{N}$ is not a direct summand of $Y$ in the coarse category,
 Corollary~\ref{cor_2-9} yields a constant $r>0$ and a sequence $(z_n)$ such that $0<d(y_n, z_n)<r$ for every $n\in \N$.

 Put $a_n =\tau^{-1}(y_n)$ and $b_n=\tau^{-1}(z_n)$ for every $n\in \N$.
 Then, $|a_n|\to \infty$ and, as in the proof of $(1)\Rightarrow(2)$,
  there exists $R>0$ such that $0< d(a_n, b_n)<R$ for every $n\in \N$.
In particular, passing to a subsequence, we may assume that $|a_{n+1}|>|a_n|+n+2R$ for every $n\in \N$.
Then the convergence $c_n\to 0$ allows us to choose a function $h\in \SO_b(X)$ such that $h(a_n)=c_n$ and $h(b_n)=0$ for every $n\in \N$
(cf. Example \ref{ex_2-1}).
Then $Th(y_n)=T\underline{c_n}(y_n)>\varepsilon_0$
 and $Th(z_n)=T\underline{0}(z_n)=0$ by Corollary~\ref{cor_2-5}.
This contradicts $Th\in \SO(Y)$.
Hence $T$ is continuous in the topology of uniform convergence.
Applying the same argument to $T^{-1}:\SO(Y)\to\SO(X)$, and using that $\mathbb{N}$ is not a direct summand of $X$ in the coarse category,
 we see that $T^{-1}$ is also continuous in the topology of uniform convergence., and hence $T$ is a homeomorphism in the topology of uniform convergence.
\end{proof}

As a consequence of this theorem, we obtain the following result,
 in which the chain-connectedness of $X$ and $Y$ in \cite[Lemma 4.4]{Iwa2025}
  is replaced by the condition that $\N$ is a direct summand of neither space in the coarse category.

\begin{corollary}\label{cor_3-02} 
Let $X$ and $Y$ be proper metric spaces such that
 $\N$ is not a direct summand of either $X$ or $Y$ in the coarse category.
If $T: \SO(X) \to \SO(Y)$ is a lattice isomorphism with
 $T\underline{0}=\underline{0}$, then the restriction of $T$ to $\SO_b(X)$ is a lattice isomorphism from $\SO_b(X)$ onto $\SO_b(Y)$.
  \hfill$\Box$
\end{corollary}

\medskip
\begin{center}
\textsc{Acknowledgement}
\end{center}
\par

This work was supported by JSPS KAKENHI Grant Number JP24K06726.

\end{document}